\documentclass[12pt]{article}
\usepackage{amssymb}
\usepackage{latexsym,bm}
\usepackage{graphicx}
\usepackage{amsmath}

\newcommand{\bea}{\begin{eqnarray*}}
\newcommand{\eea}{\end{eqnarray*}}
\newcommand{\be}{\begin{equation}}
\newcommand{\ee}{\end{equation}}
\newcommand{\ben}{\begin{eqnarray*}}
\newcommand{\een}{\end{eqnarray*}}

\date{}
\begin{document}
\title{The largest graphs with given order and diameter: A simple proof\footnote{E-mail addresses:
{\tt 235711gm@sina.com}(P.Qiao),
{\tt zhan@math.ecnu.edu.cn}(X.Zhan).}}
\author{Pu Qiao, Xingzhi Zhan\thanks{Corresponding author.}\\
{\small Department of Mathematics, East China Normal University, Shanghai 200241, China}
 } \maketitle
\begin{abstract}
A consequence of Ore's classic theorem characterizing
the maximal graphs with given order and diameter is a determination of the
largest such graphs.  We give a very short and simple proof of this smaller
result, based on a well-known elementary observation.
\end{abstract}

\vskip 6mm

We consider finite simple graphs. For terminology and notations we follow the book [2].
The {\it order} of a graph is its number of vertices, and the {\it size} its number of edges.
Denote by $V(G)$ and $E(G)$ the vertex set and edge set of a graph $G$ respectively.
We give a very short and simple proof of the following theorem of Ore [1].

{\bf Theorem} (Ore). {\it For $d\ge2$, the maximum size of a simple graph of order $n$ and
diameter $d$ is $d+(n-d-1)(n-d+4)/2$.  This size is attained by a graph  $G$ if and only if $G$
consists of a path $P$ of length $d$ such that the vertices outside $P$ form
a clique and are each adjacent to the first three or last three among some
three or four consecutive vertices on $P$.}

{\bf Proof.} Let $G$ be a simple graph of order $n$ and diameter $d.$ 
Since $G$ is of diameter $d$, there are vertices $x$ and $y$ which are at distance $d$.
Let $P$ be an $(x,y)$-path of length $d$ and denote $S=V(G)\setminus V(P)$.
To avoid bringing $x$ and $y$ closer, every vertex of $S$ has at most three neighbors on $P$, and if there
are three then they are consecutive on $P$.  Also, the $n-d-1$ vertices
of $S$ induce at most $\binom{n-d-1}2$ edges.  Counting also the edges on $P$,
we thus have $|E(G)|\le d+3(n-d-1)+\binom{n-d-1}2=d+(n-d-1)(n-d+4)/2$.

To achieve equality and thus prove sharpness of the bound, $S$ must be a
clique, and each vertex of $S$ must have three consecutive neighbors along $P$.
Since the vertices of $S$ are pairwise adjacent, their neighborhoods on $P$
together can include only at most four consecutive vertices without providing a shorter
$(x,y)$-path.  Hence the extremal graphs are formed by choosing three or four
consecutive vertices along $P$ and making each vertex of $S$ adjacent to the
first three or the last three of them. $\Box$

This proof makes it clear why there is the term $d$ and where the factor $n-d-1$ comes from in the expression
of the maximum size.  Ore [1] proved the result by first characterizing the maximal $n$-vertex graphs with diameter $d.$
Zhou, Xu and Liu [3] gave a different proof of the maximum size by considering the complement graph,
but they did not treat the extremal graphs.

Ore [1] also considered $k$-connected graphs. One would like to generalize the above argument to a simple proof of Ore's
extremal result for $k$-connected graphs with diameter $d$, but this does not
work.  The problem is that in applying Menger's Theorem [2, p.167] to obtain $k$
internally disjoint paths joining two vertices at distance $d$, some of the
paths may have length greater than $d$.  The simplest example is an odd cycle.
Ore was able to solve the more general problem by characterizing all the
diameter-critical graphs.

\vskip 5mm
{\bf Acknowledgement.} The authors are grateful to Professor Douglas B. West whose kind and detailed suggestions
have simplified an earlier version.
The research  was supported by Science and Technology Commission of Shanghai Municipality (STCSM) grant 13dz2260400 and
 the NSFC grant 11671148.

\end{document}